
\input amstex.tex
\documentstyle{amsppt}

\footline={\hss{\vbox to 2cm{\vfil\hbox{\rm\folio}}}\hss}
\nopagenumbers
\def\DJ{\leavevmode\setbox0=\hbox{D}\kern0pt\rlap
{\kern.04em\raise.188\ht0\hbox{-}}D}

\def\txt#1{{\textstyle{#1}}}
\baselineskip=13pt
\def\hf{{\textstyle{1\over2}}}
\def\a{\alpha}\def\b{\beta}
\def\d{{\,\roman d}}
\def\e{\varepsilon}
\def\f{\varphi}
\def\G{\Gamma}
\def\k{\kappa}
\def\s{\sigma}

\def\={\;=\;}

\def\zt{\zeta(\hf+it)}
\def\zr{\zeta(\hf+ir)}
  
\def\no{\noindent}  
\def\R{\Re{\roman e}\,} \def\I{\Im{\roman m}\,} 
\def\z{\zeta} 
 
\def\H{H_j^4({\txt{1\over2}})}  
\def\hf{{\textstyle{1\over2}}}
\def\txt#1{{\textstyle{#1}}}
\def\f{\varphi}

\font\tenmsb=msbm10
\font\sevenmsb=msbm7
\font\fivemsb=msbm5
\newfam\msbfam
\textfont\msbfam=\tenmsb
\scriptfont\msbfam=\sevenmsb
\scriptscriptfont\msbfam=\fivemsb
\def\Bbb#1{{\fam\msbfam #1}}

\def \NN {\Bbb N}

\def \RR {\Bbb R}
\def \ZZ {\Bbb Z}

\font\ff=cmr8
\def\txt#1{{\textstyle{#1}}}
\baselineskip=13pt

\font\teneufm=eufm10
\font\seveneufm=eufm7
\font\fiveeufm=eufm5
\newfam\eufmfam
\textfont\eufmfam=\teneufm
\scriptfont\eufmfam=\seveneufm
\scriptscriptfont\eufmfam=\fiveeufm
\def\mathfrak#1{{\fam\eufmfam\relax#1}}

\font\tenmsb=msbm10
\font\sevenmsb=msbm7
\font\fivemsb=msbm5
\newfam\msbfam
     \textfont\msbfam=\tenmsb
      \scriptfont\msbfam=\sevenmsb
      \scriptscriptfont\msbfam=\fivemsb
\def\Bbb#1{{\fam\msbfam #1}}

\def \NN {\Bbb N}

\def \RR {\Bbb R}
\def \ZZ {\Bbb Z}

  \def\rightheadline{{\hfil{\ff
  On the moments of Hecke series at central points II}\hfil\tenrm\folio}}

  \def\leftheadline{{\tenrm\folio\hfil{\ff
   A. Ivi\'c and M. Jutila}\hfil}}
  \def\emptyheadline{\hfil}
  \headline{\ifnum\pageno=1 \emptyheadline\else
  \ifodd\pageno \rightheadline \else \leftheadline\fi\fi}
   
\topmatter  
\title 
ON THE MOMENTS OF HECKE SERIES AT CENTRAL POINTS II
\endtitle
\author   Aleksandar Ivi\'c and Matti Jutila \endauthor
\address
Aleksandar Ivi\'c, Katedra Matematike RGF-a
Universiteta u Beogradu, \DJ u\v sina 7, 11000 Beograd,
Serbia (Yugoslavia).
\medskip
Matti Jutila, Department of Mathematics, University
of Turku, FIN-20014,  Turku, Finland.
\bigskip
\endaddress
\email {\tt 
aivic\@rgf.bg.ac.yu, jutila\@utu.fi} \endemail
\dedicatory Functiones et Approximatio XXXI (2003), 7-22
\enddedicatory
\abstract
{We prove, in standard notation from spectral theory,
 the asymptotic  formula ($B>0$)
$$\sum_{\k_j\le T}\a_j H_j(\hf) = \left({T\over\pi}\right)^2
- BT\log T + O(T(\log T)^{1/2}),
$$
by using an approximate functional equation for $H_j(\hf)$
and the Bruggeman-Kuznet-\break sov trace formula. We indicate how the error
term may be improved to $O(T(\log T)^\e)$.}
\endabstract
\endtopmatter

\head 
1. Introduction and statement of results
\endhead

The purpose of this paper is to continue the work begun
by the first author in [6]. Therein he obtained asymptotic formulas for
sums of $H_j^3(\hf)$ and $\H$, where $H_j(s)$ is the Hecke series
($s = \s + it$ will denote a complex variable)
$$
H_j(s)\;=\;\sum_{n=1}^\infty t_j(n)n^{-s}\qquad(\s > 1),\leqno(1.1)
$$
associated with the Maass wave form
$\psi_j(z)$, where $\rho_j(1)t_j(n) = \rho_j(n)$ and $
\rho_j(n)$ is the $n$-th Fourier coefficient of $\psi_j(z)$.
The function $H_j(s)$ can be continued to an entire function.
It  satisfies the functional equation
$$
H_j(s) = 2^{2s-1}\pi^{2s-2}\G(1-s+i\k_j)\G(1-s-i\k_j)
(\e_j\cosh(\pi\k_j)-\cos(\pi s))H_j(1-s),\leqno(1.2)
$$
where $\e_j \,(= \pm1)$ is the so-called parity sign of $\psi_j(z)$.
By $\,\{\lambda_j = \kappa_j^2 + {1\over4}\} \,\cup\, \{0\}\,$ we denote
the eigenvalues (discrete spectrum) of the hyperbolic Laplacian
$$
\Delta=-y^2\left({\left({\partial\over\partial x}\right)}^2 +
{\left({\partial\over\partial y}\right)}^2\right)
$$
acting over the Hilbert space composed of all
$\Gamma$-automorphic functions which are square integrable with
respect to the hyperbolic measure  ($\Gamma = \roman {PSL}(2,\ZZ)$).
For other relevant notation involving spectral theory the reader
is referred to [5], [6] or Y. Motohashi's comprehensive monograph [12].
The method used in [6] could not furnish the asymptotic formula
for sums of $H_j(\hf)$, but only the bounds
$$
T^2(\log T)^{-7/2} \ll \sum_{\k_j \le T}\a_jH_j(\hf) \ll T^2(\log T)^{1/2}
\leqno(1.3)
$$
were obtained, where as usual we set
$$
\a_j = |\rho_j(1)|^2(\cosh\pi\kappa_j)^{-1}.
$$

\smallskip
The aim of this paper is to improve (1.3) to a sharp
asymptotic formula, given by

\bigskip
THEOREM 1. {\it We have}
$$
\sum_{\k_j\le T}\a_j H_j(\hf) +
{2\over\pi}\int_0^T \,{|\zt|^2\over|\z(1+2it)|^2}\,\d t
= \left({T\over\pi}\right)^2
+ O(T(\log T)^{1/2}).\leqno(1.4)
$$

\smallskip\no

It remains yet to evaluate the weighted integral of the mean square of
$|\zt|$  in (1.4). The evaluation of this integral is given by

\bigskip
THEOREM 2. {\it There exist constants $A\,(>0)$ and $B$ which are 
effectively computable such that}
$$
\int_0^T \,{|\zt|^2\over|\z(1+2it)|^2}\,\d t =
T(A\log T + B) + O_\e(T^{{33\over35}+\e}).\leqno(1.5)
$$

\medskip
{\bf Corollary}. {\it If $A$ is the constant appearing in} (1.5), {\it then}
$$
\sum_{\k_j\le T}\a_j H_j(\hf) = \left({T\over\pi}\right)^2 
- {2A\over\pi}T\log T + O(T(\log T)^{1/2}).\leqno(1.6)
$$

\medskip \no
In (1.5) and later $\e$ denotes positive, arbitrarily small constants,
not necessarily the same ones at each occurrence. The formula (1.6)
shows that there are actually two main terms in the asymptotic formula
for the sum of $\a_j H_j(\hf)$.
Although the error term in (1.6) is probably too large by
a factor of $\sqrt{\log T}$, the method of proof of Theorem 1
does not allow any further improvement, if we use the weight function
(2.14). However, by a suitable choice of the weight function the error
terms in (1.4), (1.6) (and (1.7)) may be improved to $O(T(\log T)^\e)$. We
preferred to work directly with the Gaussian weight function (2.14)
because of its classical flavour. This already leads to (1.6) with
two main terms, which is the novelty of the paper.

It may be remarked that, with our method of proof,
 we can obtain the asymptotic formula
$$
\sum_{\k_j\le T}\a_j  = \left({T\over\pi}\right)^2 
+ O(T(\log T)^{1/2}).\leqno(1.7)
$$
This should be compared to a result
of N.V. Kuznetsov (see [12, p. 92] with $m = 1$), who had (1.7)
with the error term $O(T\log T)$, so that our result is somewhat
sharper.

In what concerns the true order of sums of
$\a_j H^k_j(\hf)$, it was conjectured in [6] that, for $k\in \NN$ fixed,
$$
\sum_{\k_j \le T}\a_j H_j^k(\hf) 
+ {2\over{\pi}}\int_{0}^T {|\zt|^{2k}\over|\z(1+2it)|^2}\d t
= T^2P_{{1\over2}(k^2-k)}(\log T)
+ O_{\e,k}(T^{1+c_k+\e}),\leqno(1.8)
$$
where $P_{{1\over2}(k^2-k)}(z)$ is a suitable polynomial of degree
${1\over2}(k^2-k)$ in $z$ whose coefficients depend on $k$,
and $0 \le c_k < 1$. We actually have $c_1 = c_2 = 0$,
and  even sharper results in these cases by (1.6) and
Y. Motohashi's result [11], respectively. Namely he proved the asymptotic 
formula ($\gamma = 0.5772157\ldots\,$ is Euler's constant)
$$
\sum_{\k_j\le T}\a_jH^2_j(\hf) = 2\pi^{-2}T^2(\log T + \gamma - \hf
- \log(2\pi)) + O(T\log^6T),
$$
while the proofs in [6], in the
cases $k = 3,4$, show that (1.8) holds with $c_3 = 1/7, c_4 = 1/3$.
We also note that the main term in Theorem 1, namely $(T/\pi)^2$,
is exactly of the form predicted by Random matrix
theory (see J.B. Conrey [1] and the work by J.B. Conrey et al.
[2]). This theory also gives the correct value of the leading coefficient of
the polynomial $P_{{1\over2}(k^2-k)}(z)$  for the cases $k = 2,3,4$, when
the asymptotic formulas for the sums in question are known.

\medskip
Our method of proof consists of using the Bruggeman-Kuznetsov trace formula
(cf. Lemma 1), coupled with a simple approximate functional
equation for $H_j(\hf)$ (of length $\asymp \k^2_j$) for Theorem 1
(cf. Lemma 2). This is proved in Section 2, which contains the
necessary lemmas. The crucial lemma is Lemma 3, which shows that,
in our case, the contribution of the Kloosterman sum part in the 
trace formula is negligible.
Theorem 1 is proved in Section 3, and Theorem 2
in Section 4. Finally in Section 5 we discuss how the error terms
in (1.4), (1.6) and (1.7) may be improved to $O(T(\log T)^\e)$.

\head 
2. The necessary lemmas
\endhead
{\bf Lemma 1.} (The first Bruggeman-Kuznetsov trace formula). 
{\it Let $f(r)$ be an even, regular function for $|\I r| \le \hf$
such that $f(r) \ll (1+|r|)^{-2-\delta}$ for some $\delta>0$. Then}
$$\eqalign{&
\sum_{j=1}^\infty \a_j t_j(m)t_j(n)f(\k_j) 
+ {1\over\pi}\int_{-\infty}^\infty {\s_{2ir}(m)\s_{2ir}(n)\over
(mn)^{ir}|\z(1+2ir)|^2}f(r)\d r\cr&
= {1\over\pi^2}\delta_{m,n}\int_{-\infty}^\infty r\tanh(\pi r)f(r)\d r
+ \sum_{\ell=1}^\infty{1\over\ell}
S(m,n;\ell)f_+\left({4\pi\sqrt{mn}\over\ell}\right),\cr}\leqno(2.1)
$$
{\it where $\delta_{m,n} = 1$ if $m = n$ and zero otherwise} ($m,n >0$), 
$\s_a(d) = \sum_{d\mid n}d^a$, $S(m,n;\ell)$ {\it is the Kloosterman sum and}
$$
f_+(x) = {2i\over\pi}\int_{-\infty}^\infty {r\over\cosh(\pi r)}J_{2ir}(x)
f(r)\d r.\leqno(2.2)
$$

The $J$-Bessel function is defined (see e.g., N.N. Lebedev
[9]) as
$$
J_\nu(z) = \sum_{k=0}^\infty {(-1)^k(z/2)^{\nu+2k}\over\G(k+1)\G(k+\nu+1)}
\qquad(|\arg z| < \pi).\leqno(2.3)
$$
The proof of Lemma 1 is to be found e.g. in Y. Motohashi 
[12, Chapter 2].

\medskip
{\bf Lemma 2.} {\it Let $\k_j = (1+o(1))K,\, r = (1+ o(1))K$ } ($r\in \RR$)
{\it as
$K\to\infty,  Y = (1+\delta){K^2\over4\pi^2}$, with $\delta>0$
a given constant. Then, for any fixed positive constant $A>0$, there
exists a constant $C = C(A,\delta) > 0$ such that, for $h = C\log K$,
we have}
$$
H_j(\hf) = \sum_{n\le (1+\delta)Y}t_j(n)n^{-1/2}{\roman e}^{-(n/Y)^h} 
+ O(K^{-A}),\leqno(2.4)
$$
{\it and} 
$$
\z(\hf+ir)\z(\hf-ir) = \sum_{n\le (1+\delta)Y}\s_{2ir}(n)n^{-{1\over2}-ir}
{\roman e}^{-(n/Y)^h}+ O(K^{-A}).\leqno(2.5)
$$

\medskip
{\bf Proof.} We start from the Mellin inversion integral   (see e.g.,
[4, (A.7)])
$$
{\roman e}^{-(n/Y)^h} = {1\over2\pi i}\int_{(c)}\left({Y\over n}\right)^w
\G(1 + {w\over h})\,{\d w\over w}\quad(c > 0, \,Y \gg 1),\leqno(2.6)
$$
where $\int_{(c)}$ denotes integration over the line $\R w = c$. We
use (1.1) and (see [4, Chapter 1])
$$
\z(s)\z(s-a) = \sum_{n=1}^\infty \s_a(n)n^{-s}
\qquad(\s > \max(1, \,1+\R a)), \leqno(2.7)
$$
to obtain from (2.6)
$$
\sum_{n=1}^\infty t_j(n)n^{-1/2}{\roman e}^{-(n/Y)^h}
= {1\over2\pi i}\int_{(1)}H_j(\hf+w)\G(1 + {w\over h})\,{Y^w\over w}\d w
\leqno(2.8)
$$
and
$$
\sum_{n=1}^\infty 
\s_{2ir}(n)n^{-{1\over2}-ir}{\roman e}^{-(n/Y)^h}
= {1\over2\pi i}\int_{(1)}\z(w+\hf+ir)\z(w+\hf-ir)
\G(1 + {w\over h})\,{Y^w\over w}\d w.\leqno(2.9)
$$
We shall give only the detailed proof of the more complicated formula (2.4).
The proof of (2.5) is analogous, being based on the use of (2.9). 
The series in (2.8) can be truncated at $n = (1+\delta)Y$ with the error 
$\ll K^{-A}$. On the right-hand
side of (2.8) we replace the line of integration by ${\Cal L}
= \gamma_1 \cup \gamma_2 \cup\gamma_3 \cup\gamma_4 \cup\gamma_5 $, where
$\gamma_1$ is the line from $-1-i\infty$ to $-1-ih^2$,
$\gamma_2$ is the line segment from   $-1-ih^2$ to $-\hf h - ih^2$,
$\gamma_3$ is the line segment from   $-\hf h - ih^2$ to $-\hf h + ih^2$,
$\gamma_4$ is the line segment from   $-\hf h + ih^2$ to $-1+ih^2$, and
$\gamma_5$ is the line from  $-1+ih^2$ to $-1 +i\infty$. In doing
this we pass the pole $w=0$ which, by the residue theorem, gives us
the desired contribution $H_j(\hf)$. By the functional
equation (1.2) we have
$$
H_j(\hf + w) = X_j(\hf + w)H_j(\hf - w) \leqno(2.10)
$$
with
$$
X_j(\hf + w) = (2\pi)^{2w}\pi^{-1}\G(\hf-w+i\k_j)\G(\hf-w-i\k_j)
(\e_j\cosh(\pi\k_j) + \sin(\pi w)).\leqno(2.11)
$$
To bound the gamma factors on $\Cal L$ we use Stirling's formula in
the form
$$
\G(\s + it) \ll |t|^{\s-{1\over2}}{\roman e}^{-\pi|t|/2}
\qquad(|t| \ge t_0),\leqno(2.12)
$$
which is valid uniformly for $0 \le \s \le |t|^{2/3}$.
To see this, note that
$$\eqalign{&
\R\left\{\log\G(\s+it) - \log\G(it)\right\}\cr&
= \R\left(\int_0^\s{\G'(x+it)\over\G(x+it)}\d x\right)\cr&
= \Re\left\{\int_0^\s \left(\log(x+it) - {1\over2(x+it)}
+ O\bigl({1\over(x+it)^2}\bigr)\right)\d x\right\}\cr&
\le \hf\s\log(t^2+\s^2) + O(\s t^{-2}) \le \s\log|t| +
 O((\s+\s^3)t^{-2}),\cr}
$$
hence (2.12) follows from Stirling's formula for $\G(it)$, and can be used
to bound  the gamma-factors appearing in the expression for
$X_j(\hf+w)$.

\smallskip
We have first
$$
\int_{\gamma_1}H_j(\hf+w)\G(1 + {w\over h})\,{Y^w\over w}\,\d w
\ll \int_{h^2}^\infty \exp\Bigl(-{\pi v\over2h}\Bigr)
(K^2+v^2)\d v \ll K^{-A},
$$
if $C$ in the formulation of the lemma is sufficiently large, and an 
analogous bound holds for the integral over $\gamma_5$.

\smallskip
Next, on $\gamma_2$ and on $\gamma_4$, the integrand is
$$
\ll (\k_j^2-h^4)^{-\s}(4\pi^2Y)^\s{\roman e}^{-\pi h/2} 
\ll {\roman e}^{-\pi h/2} \ll K^{-A},
$$
so that the corresponding integrals are of the desired order of magnitude.

\smallskip
Finally, on $\gamma_3$, the integrand is
$$
\ll \k_j^h(4\pi^2Y)^{-h/2} \le \left({(1+o(1))K^2\over4\pi^2Y}\right)^{h/2}
\le (1+\hf\delta)^{-h/2} \le K^{-A}
$$
for any fixed $A>0$. Combining the above bounds we obtain (2.4).

\medskip
{\bf Lemma 3.} {\it For $C\sqrt{\log K} \le G \le K$ and a
sufficiently large constant $C>0$ we have}
$$
\sum_{K\le \k_j\le K+G}\a_j H_j(\hf) \ll GK.\leqno(2.13)
$$

\medskip {\bf Proof.} First we remark that the slightly weaker bound 
$GK\sqrt{\log K}$ for the sum in (2.13) follows by applying the 
Cauchy-Schwarz inequality and the bound for sums of $\a_j$ and $\a_j
H_j^2(\hf)$ in short intervals; such bounds are given by Y. Motohashi
[12, pp. 121-122 and (3.5.13)]. 

Secondly, in the proof of Lemma 3 we may restrict $G$ to 
$G = G_0 = C\sqrt{\log K}.$
For larger $G$ we divide $[K,\,K+G]$ into $\ll G/G_0$ subintervals of
length $G_0$, to each of which we apply (2.13) with suitable $K$ and
$G = G_0$. Adding up all the results we arrive at (2.13).

The idea of proof of (2.13) is actually the same as the one that 
will be used in the proof of Theorem 1, and for the proof of Theorem
1 we need (2.13) only with $G = C\sqrt{\log K_0}, K_0 \le K \le 2K_0$. 
Lemma 3 is in fact a local  version of Theorem 1. Thus let, 
for $G = C\sqrt{\log K}$,
$$
f(r,K) 
:= {(r^2 + {\txt{1\over4}})\over(r^2 + 1000)}
\left\{\exp\left(-\left({r-K\over G}\right)^2\right)
+ \exp\left(-\left({r+K\over G}\right)^2\right)\right\}.\leqno(2.14)
$$
This function, which is a Gaussian weight function and a slightly
modified function of the function used systematically by Y. Motohashi [11],
[12], clearly satisfies the conditions of Lemma 1. 
To begin the proof, we apply Lemma 1 (taking $n=1$), combined with Lemma 2,
where $\delta>0$  is a small constant. This yields, since $H_j(\hf)\ge0$
(see S. Katok--P. Sarnak [8] for a proof),
$$
\eqalign{& \sum_{K\le \k_j\le K+G}\a_j H_j(\hf) \le
2\sum_{j=1}^\infty\a_jH_j(\hf)f(\k_j,K) \cr&
= {2\over\pi^2}\int_{-\infty}^\infty r\tanh(\pi r)f(r,K)\d r
- {2\over\pi}\int_{-\infty}^\infty
\,{|\zr|^2\over|\z(1+2ir)|^2}\,f(r,K)\d r 
\cr&
+ 2\sum_{m\le(1+\delta)^2K^2/(4\pi^2)}m^{-1/2}{\roman e}^{-(m/Y)^h}
\sum_{\ell=1}^\infty {1\over\ell}S(m,1;\ell)
f_+\left({4\pi\over\ell}\sqrt{m}\right)+ o(1)\cr&
\le {2\over\pi^2}\int_{-\infty}^\infty r\tanh(\pi r)f(r,K)\d r\cr&
+ 2\sum_{m\le(1+\delta)^2K^2/(4\pi^2)}m^{-1/2}{\roman e}^{-(m/Y)^h}
\sum_{\ell=1}^\infty {1\over\ell}S(m,1;\ell)
f_+\left({4\pi\over\ell}\sqrt{m}\right)+ o(1),\cr}\leqno(2.15)
$$
where $f_+$ is given by (2.2) with $f(r) = f(r,K)$. 

We have first
$$
\int_{-\infty}^\infty r\tanh(\pi r)f(r,K)\d r
\ll K\int_{K-G\log^2K}^{K+G\log^2K}{\roman e}^{-(r-K)^2/G^2}\d r + 1
\ll GK.\leqno(2.16)
$$

The crucial step in the proof is to show that, for any fixed $A>0$,
$$
\sum_{\ell=1}^\infty {1\over\ell}S(m,1;\ell)
f_+\left({4\pi\over\ell}\sqrt{m}\right) \ll K^{-A},\leqno(2.17)
$$
provided that we choose $G \ge C\sqrt{\log K}$.

To begin with, we may truncate the $\ell$-sum in (2.17) to the range
$1 \le \ell \le K^B$ for some constant $B>1$. To see this, we move the
line of integration in the integral defining $f_+$ (cf. (2.2)) to
$\I r = -1$. Since $f(-\hf i,K) = 0$, there is no pole of the integrand.
Then we use the series representation (see (2.3))
$$
J_{2+ix}(z) = \sum_{k=0}^\infty {(-1)^k(z/2)^{2+ix+2k}
\over\G(k+1)\G(k+ 2+ix+1)} \quad (z = 4\pi\sqrt{m}/\ell \ll K^{1-B}),
$$
which shows that the contribution of $\ell > K^{B}$ is $\ll K^{-A}$
for any fixed $A>0$, provided that $B = B(A)$ is sufficiently large.

In the remaining sum, we substitute (see e.g., [9, p. 139])
$$
J_{2ir}(x) - J_{-2ir}(x) = {2i\over\pi}\sinh(\pi r)
\int_{-\infty}^\infty \cos(x\cosh u)\cos(2ru)\d u.
$$
Integration by parts shows that, for $x>0$ and $r\ge0$,
$$
\eqalign{
J_{2ir}(x) - J_{-2ir}(x)
& = {2i\over\pi}\sinh(\pi r)
\int_{-\log^2K}^{\log^2K} \cos(x\cosh u)\cos(2ru)\d u\cr&
+ O\left(x^{-1}(r+1)\exp(\pi r - \hf\log^2K)\right).\cr}\leqno(2.18)
$$
The error term in (2.18) clearly contributes $\ll K^{-A}$ to the
sum in (2.17). The main term in (2.18) will contribute to $f_+$
$$
-{4\over\pi^2}\int_{-\log^2K}^{\log^2K} \cos(x\cosh u)\int_0^\infty
rf(r,K)\tanh(\pi r)\cos(2ru)\d r\d u.\leqno(2.19)
$$
In the inner integral we use
$$
r\tanh(\pi r) = r\,{\roman {sign}}\, r + O(|r|\exp(-\pi |r|)),\leqno(2.20)
$$
and make the change of variable $r = K + Gx$. The $x$ integral can be
truncated at $|x| = \log^2K$ with error $\ll K^{-A}$. The rational function
in $x$ in the integrand is expanded by Taylor's series, taking so many
terms that the  error  will again make
a contribution which will be $\ll K^{-A}$. Then (2.19) will become
$$
= \R \int_{-\log^2K}^{\log^2K} P(u,K,G)\cos(x\cosh u)\exp(-(G^2u^2+2iKu))\d u
+ O(K^{-A}),\leqno(2.21)
$$
where $P(u,K,G)$ is a polynomial in $u,K$ and $G$. Here we used the
familiar integral
$$
\int_{-\infty}^\infty \exp(Ax - Bx^2)\,\d x \;=\; 
\sqrt{\pi\over B}\exp\left({A^2\over4B}\right)\qquad(\R B > 0),\leqno(2.22)
$$
and $P(u,K,G)$ may be evaluated by successive differentiation of (2.22)
as the function of $A$.

If $G \ge C\sqrt{\log K}$ with large $C>0$, then the integration in (2.21)
can be restricted to the interval $|u| \le u_0$, where $u_0$ is
a small positive constant, and the error thus made will be $\ll K^{-A}$.
Then the relevant exponential factor will be of the form
$$
\exp(ig(u)),\;g(u) = \pm x\cosh u + 2Ku,\; g'(u) = \pm x\sinh u + 2K \gg K
$$
for $|x| \le BK$ and any constant $B > 0$ and $|u| \le u_0$ with sufficiently
small $u_0$, since $\sinh u = u + O(|u|^3)$ for small $u$.
In our case $x = 4\pi\sqrt{m}/\ell \le 2(1+\delta)K$ by (2.4). 
Thus the corresponding integral will have no saddle points, and by a large
number of successive integrations by parts it transpires that the integral
in question will be $\ll K^{-A}$, and so will also be
$f_+(4\pi\sqrt{m}/\ell)$. Therefore (2.17) holds, and Lemma 3 follows
from (2.15)--(2.17).

\medskip
{\bf Lemma 4.} {\it If $A(s) = \sum_{m\le M}a(m)m^{-s}$ with
$a(m) \ll_\e m^\e$, then we have }
$$\eqalign{&
\int_0^T|\zt|^2|A(\hf+it)|^2\d t\cr&
= T\sum_{h,k\le M}{a(h){\overline {a(k)}}\over hk}(h,k)
\left(\log{T(h,k)\over2\pi hk} + 2\gamma - 1\right) + E(T,A),\cr}
\leqno(2.23)
$$
{\it with $E(T,A) \ll_\e T^{1/3+\e}M^{4/3}$ if $M \ll T^C$ for some $C>0$.}

\medskip This mean value result was proved by Y. Motohashi [10].                                                                           
\head 
3. The proof of Theorem 1
\endhead
\medskip
As in the proof of Lemma 2, we let $f(r,K)$ be defined by (2.14). We
suppose additionally that $K_0 \le K \le 2K_0,$ and that $G = G(K_0)$ is
a function of $K_0$ (later we shall choose $G = C\sqrt{\log K_0}$).
We apply Lemma 1 and Lemma 2, similarly as in (2.15). Then we divide
by $\sqrt{\pi}G$ and integrate the resulting expression over $K$ from
$K_0$ to $2K_0$. It follows that
$$
\eqalign{&
\sum_{j=1}^\infty\a_jH_j(\hf)w(\k_j) + {1\over\pi}\int_{-\infty}^\infty
\,{|\zr|^2\over|\z(1+2ir)|^2}\,w(r)\d r \cr&
= {1\over\pi^2}\int_{-\infty}^\infty 
r\tanh(\pi r)w(r)\d r
+ o(1)\cr&
+ {1\over\sqrt{\pi}G}\int_{K_0}^{2K_0}\sum_{m\le(1+\delta)^2
K^2/(4\pi^2)}m^{-1/2}{\roman e}^{-(m/Y)^h}
\sum_{\ell=1}^\infty {1\over\ell}S(m,1;\ell)
f_+\left({4\pi\over\ell}\sqrt{m}\right)\d K,\cr}\leqno(3.1)
$$
where we set
$$
w(r) \;:=\; {1\over\sqrt{\pi}G}\int_{K_0}^{2K_0}f(r,K)\d K.\leqno(3.2)
$$
Since $w(r)$ is even, it suffices to consider $r\ge 0$.
From (2.14) we obtain, with the change of variable $K = r + Gx$,
$$
w(r) \;=\; {1\over\sqrt{\pi}}\int_{(K_0-r)/G}^{(2K_0-r)/G}
{\roman e}^{-x^2}\d x + O(K_0^{-2}).\leqno(3.3)
$$
If $r \in [K_0+CG\sqrt{\log K_0} ,\,2K_0-CG\sqrt{\log K_0}]$ with large $C>0$,
then the integral in (3.3) equals $1 + O(K_0^{-2})$.  If $r >
2K_0+CG\sqrt{\log K_0}$ or $r < K_0-CG\sqrt{\log K_0}$, the integral is
$O(K_0^{-2})$. Otherwise note that, for $x\ge 0$, we have $2{\roman e}^{x}
\ge 2 + 2x + x^2$, which implies that
$$
{\roman e}^{-x} \;\le\; 2(x + 1)^{-2}\qquad(x \ge 0).\leqno(3.4)
$$
Hence using (3.2)-(3.4) we obtain ($\chi_{\Cal I}(x)$ is the characteristic  
function of the set $\Cal I$), for $r \ge 0$,
$$
w(r) = \chi_{[K_0,2K_0]}(r) + O(K_0^{-2})
+ O\Bigl\{G^3(G + \min(|r - K_0|,\,|r - 2K_0|))^{-3}\Bigr\}.\leqno(3.5)
$$
\nobreak
Using (3.5) and Lemma 2 we have, for $C>0$ sufficiently large,
\goodbreak
$$
\eqalign{&
\sum_{j=1}^\infty \a_jH_j(\hf)w(\k_j) 
= \sum_{K_0-CG\sqrt{\log K_0}\le\k_j\le2K_0+CG\sqrt{\log K_0}}
\a_jH_j(\hf)w(\k_j)  + O(1)\cr&
= \sum_{K_0\le\k_j\le2K_0}\a_jH_j(\hf) + O(1)\cr&
+ O\left(G^3\sum_{{K_0-CG\sqrt{\log K_0}}\le\k_j\le K_0}\a_jH_j 
(G + K_0 - \k_j)^{-3}\right)\cr&
+ O\left(G^3\sum_{2K_0<\k_j\le 2K_0+CG\sqrt{\log K_0}}\a_jH_j 
(G +  \k_j - 2K_0)^{-3}\right)\cr&
=  \sum_{K_0\le\k_j\le2K_0}\a_jH_j(\hf)  + O(GK_0).\cr}\leqno(3.6)
$$
Similarly we obtain, since $w(r) = w(-r)$,
$$
{1\over\pi}\int_{-\infty}^\infty\,{|\zr|^2\over|\z(1+2ir)|^2}w(r)\d r
= {2\over\pi}\int_{K_0}^{2K_0}\,{|\zr|^2\over|\z(1+2ir)|^2}w(r)\d r 
+ O(GK_0),\leqno(3.7)
$$
on using $1/\z(1+it) \ll \log t$ and $\zt \ll t^{1/6}$.
Finally we have, since (2.20) holds,
$$\eqalign{
{1\over\pi^2}\int_{-\infty}^\infty 
r\tanh(\pi r)w(r)\d r &= {2\over\pi^2}\int_{K_0}^{2K_0}r\d r + O(K_0G)\cr&
= 
{1\over\pi^2}\left\{(2K_0)^2 - K_0^2\right\} + O(GK_0).\cr}\leqno(3.8)
$$
We note that the contribution of the Kloosterman-sum part in (3.1),
analogously to (2.17), is $\ll K_0^{-A}$ for any fixed $A>0$.
Therefore from (3.1) and (3.6)--(3.8) it follows that
$$
\eqalign{&
\sum_{K_0<\k_j\le2K_0}\a_jH_j(\hf) + {2\over\pi}\int_{K_0}^{2K_0}
\,{|\zr|^2\over|\z(1+2ir)|^2}\,\d r\cr&
= {1\over\pi^2}\left\{(2K_0)^2 - K_0^2\right\} + O(GK_0).\cr}\leqno(3.9)
$$
Theorem 1 follows now from (3.9) if we choose $G = C\sqrt{\log K_0}$
with a sufficiently large constant $C>0$, replace $K_0$ by $T2^{-j}$ 
and then sum over $j = 1,2,\ldots\;$. The asymptotic formula (1.7)
follows similarly as the proof of Theorem 1, if one uses the 
technique of proof of Theorem 2. One simply takes $m = n = 1$ 
in Lemma 1 and proceeds as in the
proof of Theorem 1, only the argument is simpler and the details are
thus omitted. Namely the integral in (1.4) will appear without 
$|\zt|^2$, and will be asymptotic to $CT$.

\bigskip
\head
4. The proof of Theorem 2 
\endhead
In the general problem of evaluating $\sum_{\k_j\le T}\a_jH^k_j(\hf)$
one encounters the integrals (see (1.8))
$$
I_k(T) \;:=\;\int_{0}^T {|\zt|^{2k}\over|\z(1+2it)|^2}\d t
\qquad(k\in\NN),\leqno(4.1)
$$
where $k$ is fixed. By general convexity results for Dirichlet series
one has (see K. Ramachandra [13])
$$
I_k(T) \gg_k T(\log T)^{k^2}.\leqno(4.2)
$$
Although one expects the lower bound in (4.2) to be of the correct order
of magnitude this, like in the case of the integral without the zeta-factor
in the denominator, seems at present impossible to prove for $k \ge 3$.
In fact, even for $k = 2$, when precise results on $\int_0^T|\zt|^4\d t$
are known (see e.g., [5] and [12]), an upper bound for $I_2(T)$
corresponding to the lower bound in (4.2) seems difficult to obtain and
represents an open problem.  A slightly weaker bound, namely
$I_2(T) \ll T(\log T)^4(\log\log T)^2$, follows from [14, eqs. (3.34)-(3.36)]
by a method similar to the one used in the proof of Theorem 2.

\smallskip
What we can obtain, though, is the asymptotic formula (1.5) of Theorem
2, which will be proved now. We remark that the exponent of the error
term is by no means best possible, and the use of optimal known zero-density
estimates would certainly lead to small improvements.

We start from
$$
J_1(T) := \int_{T}^{2T}\,{|\zt|^2\over|\z(1+2it)|^2}\d t
= \int_{{\Cal A}(T)} + \int_{{\Cal B}(T)}.\leqno(4.3)
$$
Here ${\Cal A}(T)$ is the subset of points $t\in [T,2T]$ such that there
are no zeros $\rho = \b + i\gamma$ of $\z(s)$ satisfying
${3\over4} \le \b \le 1,\, 2t - \log^4T\le \gamma \le 2t + \log^4T$,
and ${\Cal B}(T) = [T,2T] \,\backslash\, {\Cal A}(T)$. From M.N. Huxley's
zero-density estimate (see [4, Chapter 11])
$$
N(\s,T) = \sum_{\b\ge\s,|\gamma|\le T}1 \ll
T^{(3-3\s)/(3\s-1)}\log^CT\qquad(C>0,\,{\txt{3\over4}} \le \s \le 1)
$$
it follows that 
$$
\mu({\Cal B}(T)) \ll T^{3/5}\log^CT,\leqno(4.4)
$$
where $\mu(\cdot)$ denotes measure. Thus, by the Cauchy-Schwarz inequality
for integrals,
$$\eqalign{
\int_{{\Cal B}(T)}{|\zt|^2\over|\z(1+2it)|^2}\d t&
\le {\left\{
\int_T^{2T}{|\zt|^4\over|\z(1+2it)|^4}\d t\cdot\mu({\Cal B}(T))
\right\}}^{1/2}\cr&
\ll T^{4/5}\log^CT,\cr}
$$
where $C$ denotes generic positive constants, and where the integral
with the fourth moment of $|\zt|$ was estimated trivially as $\ll
T\log^6T$, using $1/\z(1+2it) \ll \log t$.
If $t \in {\Cal A}(T)$, then $1/\z(\s + 2it + iv) \ll_\e t^\e$ for
$\s > 3/4$ and
$|v| \le \hf \log^4 T$ (e.g., by the technique of [15, Chapter 14]). Hence
from (2.6) we obtain ($h = \log^2 T,$ $T^\e \ll Y \ll T^{1/2}$)
$$
\eqalign{&
\sum_{n=1}^\infty\mu(n)n^{-1-2it}{\roman e}^{-(n/Y)^h} 
= {1\over2\pi i}\int\limits_{(1)}{Y^w\over\z(1+2it+w)}\G(1+{w\over h})
{\d w\over w}\cr&
= {1\over2\pi i}\int\limits_{\R w = 1,|\I w|\le \hf h^2}
{Y^w\over\z(1+2it+w)}\G(1+{w\over h}){\d w\over w}  + O(T^{-10})\cr&
= {1\over\z(1+2it)} + {1\over2\pi i}\int\limits_{\R w = 
\e-{1\over4},|\I w|\le \hf h^2}
{Y^w\over\z(1+2it+w)}\G(1+{w\over h}){\d w\over w}  + O(T^{-10})\cr&
= {1\over\z(1+2it)} + O(Y^{-1/4}T^\e)   + O(T^{-10}).\cr}\leqno(4.5)
$$
Set $a(m) = \mu(n)$ if $m = n^2$ and $a(m) = 0$ otherwise. From (4.5)
it follows that, for $t \in {\Cal A}(T)$,
$$
{1\over\z(1+2it)} = \sum_{m\le4Y^2}a(m)m^{-1/2-it}\exp(-(\sqrt{m}/Y)^h)
+ O(T^\e Y^{-1/4}).\leqno(4.6)
$$
We then obtain, using (4.4), (4.6) and the Cauchy-Schwarz inequality,
$$
\eqalign{
\int_{{\Cal A}(T)}\ldots\d t
&= \int_T^{2T}|\zt|^2\Bigl|\sum_{m\le4Y^2}a(m)
m^{-1/2-it}\exp(-(\sqrt{m}/Y)^h)\Bigr|^2\d t
\cr&
+ O_\e(T^{1+\e}Y^{-1/4}) + O(T^{4/5}\log^CT).
\cr}
$$
To evaluate the last integral we use (2.23) of Lemma 4. We obtain
$$
\eqalign{&
\int_0^T|\zt|^2
\Bigl|\sum_{m\le4Y^2}a(m)
m^{-1/2-it}\exp(-(\sqrt{m}/Y)^h)\Bigr|^2\d t \cr&
= T\sum_{\ell,k\le2Y}{\mu(\ell)\mu(k)\over \ell^2k^2}
{\roman e}^{-(\ell/Y)^h-(k/Y)^h}(\ell,k)^2\left(\log{T(\ell,k)^2\over
2\pi \ell^2k^2} + 2\gamma -1 \right)\cr&
+ O_\e(T^{1/3+\e}Y^{8/3}).\cr}
$$
Setting $d = (\ell,k), \ell = d\ell_1, k = dk_1, (\ell_1,\,k_1) = 1$,
we see that the double sum above equals
$$\eqalign{
\sum_{d\le 2Y}{\mu^2(d)\over d^2}\sum_{k_1\le{2Y\over d},
\ell_1\le{2Y\over d},(k_1,\ell_1)=(k_1,d)=(\ell_1,d)=1}
{\mu(k_1)\mu(\ell_1)\over k_1^2\ell_1^2}\times\cr
\times {\roman e}^{-(d\ell_1/Y)^h-(dk_1/Y)^h}
\left\{\log\left({T\over2\pi k_1^2\ell_1^2d^2}\right)+2\gamma-1\right\}. \cr}
$$
The terms $k_1 > Y/(2d)$, and then $\, \ell_1 > Y/(2d)$ 
are estimated trivially, 
producing an error which is $O(TY^{-1}\log^2T)$. In the remaining terms we
get rid of the exponential factor by using ${\roman e}^{-x} = 1 + O(x)$
for $x > 0$. In the inner sum we extend the summation to all 
$k_1,\,\ell_1,\,$ obtaining again an error which is $O(TY^{-1}\log^2T)$,
and similarly we extend the summation over all $d$. Finally we obtain
that the double sum above equals
 $$
 A\log T + B + O\left({\log^2T\over Y}\right)\qquad(A>0),
 $$
 where the constants $A$ and $B$ may be explicitly evaluated. Putting
 together all the expressions we wind up with
 $$
 \eqalign{
 \int_0^T\,{|\zt|^2\over|\z(1+2it)|^2}\,\d t& = T(A\log T + B)\cr&
 + O_\e(T^{1/3+\e}Y^{8/3}) + O_\e(T^{1+\e}Y^{-1/4}) + O(T^{4/5}\log^CT).\cr}
 $$
 The choice $Y = T^{8/35}$ completes the proof of (1.5) of Theorem 2.

\bigskip
\head
5. The choice of the weight function 
\endhead
We shall discuss now how the error terms in (1.4) (and thus also in (1.6)
and (1.7)) can be improved to $O(T(\log T)^\e)$. Let $S_\alpha^\beta$
be the class of smooth functions  $f(x) \,(\in C^\infty)$ introduced  
by I.M. Gel'fand and G.E. Shilov [3]. The
functions $f(x)$ satisfy for any real $x$ the inequalities
$$
\vert x^kf^{(q)}(x)\vert  \le CA^kB^qk^{k\alpha}q^{q\beta}
\qquad (k,q = 0,1,2,\ldots) \leqno(5.1)
$$ \noindent
with suitable constants $A,B,C > 0$ depending on $f$ alone. For
$\alpha = 0$ it follows that $f(x)$ is of bounded support,
namely it vanishes for $\vert x\vert \ge A$. For $\alpha > 0$
the condition (5.1) is equivalent (see [3]) to the condition
$$
\vert f^{(q)}(x)\vert \le CB^qq^{q\beta}\exp(-a\vert x\vert^{1/\alpha})
\qquad (a = \alpha/(eA^{1/\alpha})) \leqno(5.2)
$$ \noindent
for all $x$ and $q \ge 0$. We shall denote by $E_\alpha^\beta$
the subclass of $S_\alpha^\beta$ with $\alpha > 0$ consisting of
 even functions $f(x)$ such that $f(x)$ is not the
zero-function. It is shown in [3] that $S_\alpha^\beta$
is non-empty if $\beta \ge 0$ and $\alpha + \beta \ge 1$. If
these conditions hold then $E_\alpha^\beta$ is also non-empty, since
$f(-x) \in S_\alpha^\beta$ if  $f(x) \in S_\alpha^\beta$, and
$f(x) + f(-x)$ is always even.  If
$$
\hat{f}(x) = \int_{-\infty}^\infty f(u){\roman e}^{iux}\d u
$$
denotes the Fourier transform of $f(x)$, then
a fundamental property of the class $S_\alpha^\beta$ (see op. cit.) is that
$\widehat{S_\alpha^\beta} = S_\beta^\alpha$, where in general
$\widehat{U} = \{\widehat{f}(x) : f(x) \in U\}$. Henceforth
let $\f(x)\in E_{1-\delta}^{\delta}$ be non-negative, 
where $\delta > 0$ is a small constant, and set
$$
f_\f(r) = f_\f(r,K) = {r^2+{1\over4}\over r^2+1000}
\left\{\f\left({r+K\over G}\right) + \f\left({r-K\over G}\right)\right\},
\leqno(5.3)
$$
where
$$
C(\log K)^\delta \;\le G \le \sqrt{K},\qquad(C = C(\delta)>0).\leqno(5.4)
$$
The function $\f(x)$ is of fast decay by (5.2), and moreover by the 
general theory (op. cit.) the analytic continuation of $\f(z)$
certainly exists in the strip $|y| = |\I z| \le C\;(C>0)$, 
where it is of rapid decay, so that  
$f_\f(r) $ satisfies the assumptions of Lemma 1.

\medskip
Our main task is to show that (2.17) holds with $f_+$ (cf. (2.2)) relating
to $f_\f(r)$, as given by (5.3), and $G$ satisfying (5.4), where
of course it is the lower bound that is critical. We follow the reasoning
given from (2.18)--(2.22) in the proof of Lemma 3, but make the
following observations. The reason $G = C\sqrt{\log K}$ was the limit
in Lemma 3 (and indirectly in the proof of Theorem 1) is the appearance
of $\exp(-(G^2u^2 + 2iKu))$ in (2.21). 
With $f_\f(r)$ replacing $f$ (cf. (2.14)), the integral over $r$ in
(2.18) can be truncated at $|r| = \log^2 K$ with negligible error. 
While the term $2iKu$ in (2.21) (which comes
after the change of variable $r = K + Gx$) cannot be avoided, the
term $-G^2u^2$ comes  from the fact that essentially ${\roman
e}^{-x^2} \;(\in S_{1/2}^{1/2})$ is the Fourier transform of itself, which
is embodied in the formula (2.22). This factor sets the lower bound
$G = C\sqrt{\log K}$. However, in this new situation we shall obtain,
instead of $\exp(-G^2u^2)$, the function $\hat{\f}_f(x)\in S_\delta
^{1-\delta}$, which by (5.2) satisfies
$$
\hat{\f}_f(Gu) \ll \exp(-a|Gu|^{1/\delta}).\leqno(5.5)
$$
Thus we may truncate the integration in the analogue of (2.21) now
at $|u| \le u_0$, provided that $G \ge C(\log K)^\delta,  C = C(\delta)
> 0$ sufficiently large, and the analogue of (2.17) will hold again.

\medskip
It only remains to check that the integration over $[K_0, \,2K_0]$ in
the proof of Theorem 1 will go through. To do this, instead of (3.2)
consider
$$
w_\f(r) \;:=\; {1\over BG}\int_{K_0}^{2K_0} f_\f(r,K)\d K,\leqno(5.6)
$$
where $B = \hat{\f}(0) = 
\int_{-\infty}^\infty \f(x)\d x$. Since $\f(x)\in E_{1-\delta}
^\delta$, we have
$$
\f(x) \ll \exp(-a|x|^{1/(1-\delta)})\qquad(a>0).
$$
Therefore by using e.g., the inequality
$$
{\roman e}^{-x} \;\le\;24(x + 1)^{-4}\qquad(x \ge 0),
$$
we obtain the analogue of (3.5) for $w_\f(r)$. This means that the choice
$G = C(\log K_0)^\delta$ is permissible in the proof of Theorem 1, which
ends our discussion.

\vfill
\eject\topskip2cm
\Refs
\bigskip

\item{[1]} J.B. Conrey, $L$-functions and random matrices, in ``Mathematics
Unlimited" (Part I), B. Engquist and W. Schmid eds., Springer,
2001, pp. 331-352.

\item{[2]} J.B. Conrey, D.W. Farmer, J.P. Keating, M.O. Rubinstein
and N.C. Snaith, Integral moments of $L$-functions, preprint, 58pp,
arXiv:math.NT/0206018,

{\tt http://front.math.ucdavis.edu/mat.NT/0206018}.

\item {[3]} I.M. Gel'fand and G.E. Shilov,  Generalized functions
(vol. 2), Academic Press, New York-London, 1968.

\item{[4]} A. Ivi\'c, The Riemann zeta-function, John Wiley \&
Sons, New York, 1985.

\item{[5]} A. Ivi\'c,  The mean values of the Riemann zeta-function, Tata 
        Institute of Fundamental Research, Lecture Notes {\bf82}, 
    Bombay 1991 (distr. Springer Verlag, Berlin etc.).

\item{[6]} A. Ivi\'c, Moments of Hecke series at central points,
Functiones et Approximatio {\bf30}(2002), 49-82.

\item{[7]} M. Jutila and Y. Motohashi, A note on the mean value
of the zeta and $L$-functions XI, Proc. Japan Acad. {\bf78}, Ser. A
(2002), 1-6.

\item {[8]} S. Katok and P. Sarnak, Heegner points, cycles and Maass
forms,  Israel J. Math. {\bf84}(1993), 193-227.

\item{[9]} N.N. Lebedev, Special functions and their applications,
Dover,  New York, 1972.

\item{[10]} Y. Motohashi, A note on the mean value of the zeta and
$L$-functions V, Proc. Japan Acad. Ser. A {\bf52}(1986), 399-411.

\item{[11]} Y. Motohashi, Spectral mean values of Maass wave form
$L$-functions, J. Number Theory {\bf42}(1992), 258-284.

\item{[12]} Y. Motohashi, Spectral  theory of the Riemann zeta-function,
Cambridge University Press, 1997.

\item{[13]} K. Ramachandra, On the mean-value and omega-theorems for
the Riemann zeta-function, Tata Institute of Fundamental
 Research, Bombay, 1995 (distr. by Springer Verlag, Berlin etc.).

\item{[14]} K. Ramachandra and A. Sankaranarayanan, On an
asymptotic formula of Srinivasa Ramanujan, Acta Arith. {\bf109}(2003),
349-357.

\item{[15]} E.C. Titchmarsh, The theory of the Riemann zeta-function
(2nd ed.),  University Press, Oxford, 1986.

\bigskip

Aleksandar Ivi\'c

Katedra Matematike RGF-a

Universitet u Beogradu

\DJ u\v sina 7, 11000 Beograd, Serbia

{\tt aivic\@rgf.bg.ac.yu}

\bigskip

Matti Jutila

Department of Mathematics

 University of Turku, FIN-20014

  Turku, Finland.

 {\tt jutila\@utu.fi}

\endRefs
\vfill


\bye